\documentclass[11pt]{article}

\usepackage{amsmath,amssymb}
\usepackage{color}

\newtheorem{theorem}{Theorem}
  \newtheorem{lemma}{Lemma}
  \newtheorem{proposition}{Proposition}
  
  \newtheorem{definition}{Definition}
  \newtheorem{remark}{Remark}

  \newcommand{\pf}{{\bf Proof: }}

\def\Z{{\mathbb{Z}}}

\def\N{{\mathbb{N}}}
\def\Q{{\mathbb{Q}}}

\def\be{\begin{eqnarray*}}
\def\ee{\end{eqnarray*}}
\def\bee{\begin{eqnarray}}
\def\eee{\end{eqnarray}}

\begin{document}

\title{\Large{\bf A Pellian equation with primes and applications to $D(-1)$-quadruples}}
{\let\thefootnote\relax\footnotetext{Authors were supported by the Croatian Science Foundation under the project no. 6422.
A.D. acknowledges support from the QuantiXLie Center of Excellence.}}

\author{Andrej Dujella, Mirela Juki\'c Bokun and Ivan Soldo}

\date{}
\maketitle

\begin{abstract}
\noindent
In this paper, we prove that the equation $x^2-(p^{2k+2}+1)y^2=-p^{2l+1}$, $l \in \{0,1,\dots,k\}, k \geq 0$, where $p$ is an odd prime number,
is not solvable in positive integers $x$ and $y$. By combining that result with other known
results on the existence of Diophantine quadruples, we are able to prove results on the extensibility of some $D(-1)$-pairs to
quadruples in the ring $\Z[\sqrt{-t}], t>0$.
\end{abstract}

\medskip
{\small

}

\section{Introduction}

Diophantus of Alexandria raised the problem of finding four positive rational numbers $a_1, a_2, a_3, a_4$
such that $a_i a_j+1$ is a square of a rational number for each $i,j$ with $1\leq i< j\leq4$ and gave a solution
$\{\frac{1}{16}, \frac{33}{16}, \frac{17}{4}, \frac{105}{16} \}$. The first example of such a set in the ring of integers
was found by Fermat and it was the set $\{1, 3, 8, 120\}$. Replacing ``$+1$'' by ``$+n$'' suggests the following general definition:
\begin{definition}\label{def_1}
Let $n$ be a non-zero element of a commutative ring $R$. A {\it Diophantine m-tuple
with the property} $D(n)$, or simply {\it a} $D(n)$-{\it m-tuple},
is a set of $m$ non-zero elements of $R$ such that if $a, b$ are any two distinct elements from this
set, then $ab + n = k^2$, for some element $k$ in $R$.
\end{definition}

Let $p$ be an odd prime and $k$ a non-negative integer. We consider the Pellian equation
\bee\label{eq1}
x^2-(p^{2k+2}+1)y^2=-p^{2l+1}, \quad l \in \{0,1,\dots,k\}.
\eee

The existence of positive solutions of the above equation is closely related to the existence of a Diophantine quadruple in certain ring.
More precisely, the entries in a Diophantine quadruple are strictly restricted in that they appear as coefficients of
three generalized Pell equations that must have at least one common solution in positive integers.



According to Definition~\ref{def_1}, we will look at the case $n = -1$. Research on $D(-1)$-quadruples is quite active. It is conjectured that $D(-1)$-quadruples do not exist in integers (see \cite{DujeFib}).
Dujella, Filipin and Fuchs in \cite{DFF} proved that there are at most finitely many $D(-1)$-quadruples, by giving an upper bound of $10^{903}$
for their number. There is a vast literature on improving that bound (e.g., see \cite{FiFu1, BoCiMm,EsFiFu}). Very recently, in \cite{Trudg}
Trudgian proved that there are at most $3.01\cdot 10^{60}$ $D(-1)$-quadruples. In \cite{BoBuCiMm}, it is announced that the bound can be reduced to $2.5 \cdot 10^{60}$.

Concerning the imaginary quadratic fields, Dujella (see \cite{DujeGlas}) and Franu\v{s}i\'c (see \cite{Fran4})
considered the problem of existence of $D(-1)$-quadruples in Gaussian integers. Moreover, in \cite{FranKreso} Franu\v{s}i\'c
and Kreso showed that the Diophantine pair $\{1,3\}$ cannot be extended to a Diophantine quintuple in the ring $\Z[\sqrt{-2}]$.
Several authors contributed to the characterization of elements $z$ in $\Z[\sqrt{-2}]$ for which a
Diophantine quadruple with the property $D(z)$ exists (see \cite{Abu,DujeSoldo,Soldo_Misc}).  The problem of Diophantus for
integers of the quadratic field $\Q(\sqrt{-3})$ was studied in \cite{FranSoldo}. In \cite{Soldo_Studia, Soldo_Malays}, Soldo studied
the existence of $D(-1)$-quadruples of the form $\{1,b,c,d\}, b \in \{2,5,10,17,26,37,50\}$, in the ring $\Z[\sqrt{-t}], t>0$.

The aim of the present paper is to obtain results about solvability of the equation~\eqref{eq1} in positive integers.
Since we know that it is  closely related to the existence of $D(-1)$-quadruples, we use obtained results to
prove our results on the extensibility of some $D(-1)$-pairs to quadruples in the ring $\Z[\sqrt{-t}], t>0$.
of integers in the quadratic field $\Q(\sqrt{-t})$.

\section{Pellian equations}

The goal of this section is to determine all solutions in positive integers of the equation~\eqref{eq1}, which is the crucial step in proving our results in the next section. For this purpose,
we need the following result on Diophantine approximations.
\begin{theorem}[{\cite{Worley,DujeTatra}}] \label{WoDu}
Let $\alpha$ be a real number and let $a$ and $b$ be coprime non-zero integers, satisfying the inequality
\be
\left|\alpha-\frac{a}{b}\right| < \frac{c}{b^2},
\ee
where $c$ is a positive real number. Then $(a, b) = (rp_{m+1} \pm up_m, rq_{m+1} \pm uq_m)$,
for some $m \geq -1$ and non-negative integers $r$ and $u$ such that $ru < 2c$. Here $p_m/q_m$
denotes the $m$-th convergent of continued fraction expansion od $\alpha$.
\end{theorem}
If $\alpha=\frac{s+\sqrt{d}}{t}$ is a quadratic irrational, then the simple continued fraction expansion
of $\alpha$ is periodic. This expansion can be obtained by using the following algorithm.
Let $s_0 = s, t_0 = t$ and
\bee\label{c_frac}
a_n= \left\lfloor \frac{s_n+\sqrt{d}}{t_n} \right\rfloor, \quad s_{n+1}=a_n t_n - s_n, \quad t_{n+1}=\frac{d-s_{n+1}^2}{t_n}, \quad \mbox{for } n \geq 0
\eee
(see \cite[Chapter 7.7]{Niven}). If $(s_j, t_j)=(s_k, t_k)$ for $j<k$, then
\[
\alpha=[a_0, \dots, a_{j-1}, \overline{a_j, \dots, a_{k-1}}].
\]
We will combine Theorem~\ref{WoDu} with the following lemma:
\begin{lemma} [\text{\cite[Lemma 2]{DuBo}}]\label{lemaDB}
Let $\alpha$, $\beta$ be positive integers such that
  $\alpha\beta$ is not a perfect square, and
let $p_n/q_n$ denotes the $n$-th convergent of continued fraction expansion of
$\sqrt{\frac{\alpha}{\beta}}$. Let the sequences $(s_n)$ and $(t_n)$ be defined
by \eqref{c_frac} for the quadratic irrational $\frac{\sqrt{\alpha \beta}}{\beta}$. Then
\be
\alpha (r q_{n+1}+u q_n)^2 - \beta (r p_{n+1} + u p_n)^2 = (-1)^n (u^2 t_{n+1} + 2 r u s_{n+2} - r^2 t_{n+2}),
\ee
for any real numbers $r,u$.
\end{lemma}
The next lemma will be usefull, too.
\begin{lemma}[\text{\cite[Lemma 2.3.]{Fujita}}]\label{Fujita}
Let $N$ and $K$ be integers with $1<|N|\leq K$. Then the Pellian equation
\[
X^2-(K^2+1)Y^2=N
\]
has no primitive solution.
\end{lemma}
The solution $(X_0, Y_0)$ is called primitive if $\gcd(X_0, Y_0)=1$.
Now we formulate the main result of this section.
\begin{theorem}\label{tm1}
Let $p$ be an odd prime and $k$ a non-negative integer. The equation
\bee\label{eq_nova}
x^2-(p^{2k+2}+1)y^2=-p^{2l+1}, \quad l \in \{0,1,\dots, k\}
\eee
has no solutions in positive integers $x$ and $y$.
\end{theorem}

In proving  Theorem~\ref{tm1}, we will apply the following technical lemma.

\begin{lemma}\label{lema1}
If $(x,y)$ is a solution of the equation
\bee\label{eq_novak}
x^2-(p^{2k+2}+1)y^2=-p^{2k+1},
\eee
and $y\geq p^{\frac{2k+1}{2}}$, then the inequality
\[
\sqrt{p^{2k+2}+1}+\frac{x}{y}>2p^{k+1}
\]
holds.
\end{lemma}

\noindent
{\bf Proof}:
From \eqref{eq_novak} we have
\bee\label{eeq_1}
\frac{x^2}{y^2}= p^{2k+2}-\frac{p^{2k+1}}{y^2}+1.
\eee
Thus we have to consider when the inequality
\[
p^{2k+2}-\frac{p^{2k+1}}{y^2}+1 > \left(2p^{k+1} - \sqrt{p^{2k+2}+1}\right)^2
\]
is satisfied. This inequality is equivalent to 
\bee\label{eq_novakk}
\frac{p^k}{y^2}<4\left(\sqrt{p^{2k+2}+1}-p^{k+1}\right).
\eee
For $x>1$, the inequality $(1+\frac{1}{x})^\frac{1}{2} > 1 + \frac{1}{2x}-\frac{1}{8x^2}$ holds. Thus we have
\be
&&4\left(\sqrt{p^{2k+2}+1}-p^{k+1}\right)\\
&&\qquad = 4p^{k+1} \left(  \left( 1+\frac{1}{p^{2k+2}}\right)^{\frac{1}{2}}-1\right)\\
&&\qquad> 4p^{k+1} \left( \frac{1}{2p^{2k+2}} - \frac{1}{8p^{4k+4}} \right)\\
&&\qquad >  4p^{k+1} \cdot  \frac{1}{4p^{2k+2}}\\
&&\qquad =  \frac{1}{p^{k+1}}.
\ee
Since $y\geq p^{\frac{2k+1}{2}}$, i.e.
\[
\frac{1}{p^{k+1}} \geq \frac{p^k}{y^2},
\]
 we conclude that the inequality \eqref{eq_novakk} holds.
\hfill $\Box$




\medskip
\noindent
{\bf Proof of Theorem~\ref{tm1}:}

\noindent
Case 1. Let $2l+1\leq k+1$, i.e., $l \leq \frac{k}{2}$.

By Lemma~\ref{Fujita}, we know that the equation \eqref{eq_nova} has no primitive solutions. Assume that there exists a non-primitive solution $(x,y)$.
Then $p|x$ and $p|y$, so there exist $0<i\leq l, x_1, y_1 \geq 0, \gcd(x_1,y_1)=1$ such that $x=p^i x_1, y=p^i y_1$. After dividing by $p^{2i}$ in \eqref{eq_nova},
we obtain
\be\label{eq9}
x_1^2-(p^{2k+2}+1)y_1^2=-p^{2l-2i+1}, \quad 0< 2l-2i+1 \leq k+1.
\ee
But such $x_1$, $y_1$ do not exist according to Lemma~\ref{Fujita}, so we obtained a contradiction.

\noindent
Case 2. Let $2l+1= 2k+1$, i.e., $l=k$.

Let us suppose that there exists a solution $(x,y)$ of the equation \eqref{eq1}
such that $y\geq p^{\frac{2k+1}{2}}$.
Then by applying \eqref{eeq_1} we obtain
\be
|\sqrt{p^{2k+2}+1}-\frac{x}{y}|
&=&|p^{2k+2}-\frac{x^2}{y^2}+1|\cdot|\sqrt{p^{2k+2}+1}+\frac{x}{y}|^{-1}\\
&=&\frac{p^{2k+1}}{y^2}\cdot|\sqrt{p^{2k+2}+1}+\frac{x}{y}|^{-1}.
\ee
Lemma~\ref{lema1} implies
\bee\label{eqqq1}
|\sqrt{p^{2k+2}+1}-\frac{x}{y}| < \frac{p^k}{2y^2}.
\eee
Assume that $x=p^t x_1, y=p^t y_1$, where $t, x_1, y_1$ are non-negative integers and $\gcd(x_1, y_1)=~1$. Now the equation \eqref{eq1} is equivalent to
\bee\label{eqqq2}
x_1^2-(p^{2k+2}+1)y_1^2=-p^{2k-2t+1}.
\eee
Since $y\geq y_1$, from \eqref{eqqq1} we obtain
\be
|\sqrt{p^{2k+2}+1}-\frac{x_1}{y_1}| < \frac{p^k}{2y_1^2}.
\ee
Now, Theorem~\ref{WoDu} implies that
\bee \label{eqqq3}
(x_1,y_1)= (rp_{m+1} \pm up_m, rq_{m+1} \pm uq_m),
\eee
for some $m \geq -1$ and non-negative integers $r$ and $u$ such that
\bee\label{eq2}
ru < p^k.
\eee
Since $x_1$ and $y_1$ are coprime, we have $\gcd(r,u)=1$.

The terms $p_m/q_m$ are convergents of the continued fraction expansion of $\sqrt{p^{2k+2}+1}$.
Since
\[
\sqrt{p^{2k+2}+1} = [p^{k+1}, \overline{2p^{k+1}}],
\]
the period of that continued fraction expansion (and also of the corresponding sequences $(s_n)$ and $(t_n)$) is equal to 1, according to Lemma~\ref{lemaDB}, we have to consider
only the case $m=0$. We obtain
\bee\label{eq3}
(p^{2k+2}+1)(rq_1 \pm uq_0)^2-(rp_1 \pm up_0)^2 = u^2 t_1 \pm 2ru s_2 - r^2 t_2,
\eee
where
\[
s_2=p^{k+1}, \quad t_1=t_2=1, \quad p_0=p^{k+1}, p_1=2p^{2k+2}+1, \quad q_0=1, q_1=2p^{k+1}.
\]
Since the observation is similar in both signs, in what follows our focus will be to the positive sign. By comparing
\eqref{eqqq2} and \eqref{eq3}, we obtain the equation
\bee\label{eq4}
u^2-r^2 + 2ru p^{k+1} = p^{2k-2t+1}.
\eee
Now, we consider the solvability of \eqref{eq4}. 

If $r=0$, then $u^2=p^{2k-2t+1}$, and so $p$ has to be a square, which is not possible.

If $u=0$, we obtain $-r^2=p^{2k-2t+1}$, and that is not possible, too.


If $r=u$,  we have  $p^{k-2t}=2r^2$. Since $p$ is an odd prime, that is not possible.

Let $r \neq u, r,u\neq 0$. If $t<\frac{k}{2}$, then from \eqref{eq4} we conclude that $p^{k+1}|u^2-r^2$.
If $p|u+r$ and $p|u-r$, then $p|2\gcd(r,u)$, i.e., $p|2$ which is not possible.
Therefore, $p^{k+1}$ divides exactly one of the numbers $u+r$ and $u-r$. In both cases, it follows that $u+r\geq p^{k+1}$.
That implies
\[
ur\geq u+r-1>p^k,
\]
which is a contradiction with \eqref{eq2}.

Now, let us suppose that $t\geq\frac{k}{2}$. Since the equation \eqref{eq1} is equivalent to \eqref{eqqq2} and $0<2k-2t+1\leq k+1$,
by Case 1 it has no solutions.

It remains  to consider the case $y<p^{\frac{2k+1}{2}}$. Assume that there exists a solution of the equation \eqref{eq1} with this property. In that case
we can generate increasing sequence of infinitely many  solutions of the equation \eqref{eq1}. Therefore, a solution $(x,y)$ such that $y\geq p^{\frac{2k+1}{2}}$ will appear.
This contradicts with the first part of the proof of this case.

\noindent
Case 3. Let $k+1 < 2l+1 < 2k+1$, i.e., $\frac{k}{2}<l<k$.

In this case, if we suppose that the equation \eqref{eq_nova} has a solution, then multiplying that solution by $p^{k-l}$ we obtain the solution of the equation
\[
x^2-(p^{2k+2}+1)y^2=-p^{2k+1},
\]
which is not solvable by Case 2. That is the contradiction, and this completes the proof of Theorem~\ref{tm1}.
\hfill $\Box$

\begin{proposition}
Let $p=2$.
\begin{itemize}
\item[i)] If $k\equiv 0\pmod 2$, then the equation \eqref{eq_nova} has no solutions.
\item[ii)] If $k\equiv 1\pmod 2$, then in case of $l>\frac{k}{2}$ the equation \eqref{eq_nova} has a solution, and in case of $l\leq \frac{k}{2}$ it has no solutions.
\end{itemize}
\end{proposition}

\noindent
{\bf Proof}:
i) If $k\equiv 0\pmod 2$, then the equation \eqref{eq_nova} is not solvable modulo 5.

\noindent
ii) Let $k\equiv 1\pmod 2$. If $l>\frac{k}{2}$, the equation \eqref{eq_nova} has the solution of the form
\[
(x,y)=(2^{\frac{2l-k-1}{2}}(2^{k+1}-1), 2^{\frac{2l-k-1}{2}}),
\]
and therefore infinitely many solutions.

If $l\leq \frac{k}{2}$, then $2l+1\leq k+1$ and we can proceed as in Case 1 of Theorem~\ref{tm1} and conclude that the equation \eqref{eq_nova} has no solutions.
\hfill $\Box$


\section{Application to $D(-1)$-triples}

By using results from the previous section and known results on Diophantine $m$-tuples, in this section we
present the results on extensibility of certain Diophantine pairs to quadruples, in the ring  $\Z[\sqrt{-t}], t>0$.

The following result is proved in \cite{Soldo_Malays}:
\begin{theorem}[\text{\cite[Theorem 2.2]{Soldo_Malays}}]\label{tm2.2}
Let $t>0$ and $\{1,b,c\}$ be $D(-1)$-triple in the ring $\Z[\sqrt{-t}]$.
\begin{itemize}
\item[$($i$)$] If $b$ is a prime, then $c\in \Z$.
\item[$($ii$)$] If $b=2 b_1$, where $b_1$ is a prime, then $c\in \Z$.
\item[$($iii$)$] If $b=2 b_2^2$, where $b_2$ is a prime, then $c\in \Z$.
\end{itemize}
\end{theorem}

\begin{remark}\label{rem_2}
In the proof of \cite[Theorem 2.2]{Soldo_Malays}, it was shown that for every $t$ there exists such $c>0$, while the case $c<0$ is possible
only if $t|b-1$ and the equation
\bee \label{eq7}
x^2-b y^2=\frac{1-b}{t}
\eee
has an integer solution.
\end{remark}


Let $p$ be an odd prime and $b=2p^k, k\in \mathbb{N}$. We consider the extendibility of $D(-1)$-triples
of the form $\{1,b,c\}$  to quadruples in the ring $\Z[\sqrt{-t}], t>0$.
How complex that problem can be, depends on the number of divisors $t$ of $b-1$. As $b$ grows, we can
expect the larger set of $t$'s, and for each $t$ we have to consider whether there exists a solution of the
equation~\eqref{eq7}. If it is true, then the problem is reduced to solving the systems of simultaneous Pellian equations.
A variety of different methods have been used to study such kind of problems, including linear forms
in  logarithms, elliptic curves, theory around Pell's equation, elementary methods, separating the
problem into several subproblems depending on the size of parameters, etc. A sursey on that subject is given in \cite{DujeNotice}.

Therefore, since  $b-1=2p^k-1$ has to be a square, to reduce the number of $t$'s, we consider the equation of the form
\bee\label{eq_opca}
2p^k-1=q^{2j}, \quad j>0,
\eee
where $q$ is an odd prime. According to \cite[Lemma 2.9]{Koso},
if $k>1$ the equation \eqref{eq_opca} has solutions only for $(k,j)\in\{(2,1),(4,1)\}$. If $(k,j)=(2,1)$,
we obtain the Pellian equation in primes. So far known prime solutions are
$
(p,q) \in \{(5,7), (29,41),(44560482149, 63018038201),\\ (13558774610046711780701,19175002942688032928599)\}
$ (see \cite{Web_primes}).
If $(k,j)=(4,1)$, the only solution is $(p,q)=(13,239)$.

Let $k=1$. Supose that $j=mn$, where $n$ is an odd number. Then we have
\[
2p=q^{2j}+1=q^{2mn}+1=(q^{2m}+1)((q^{2m})^{n-1} -(q^{2m})^{n-2}+\dots-q^{2m}+1).
\]
Since $q$ is an odd prime, we conclude that $q^{2m}+1=2p=q^{2mn}+1$. This implies that $n=1$. This means that the only possibility for $2p=q^{2j}+1$ is that $j$ is a non-negative power of $2$.

Note that in all possible cases of $k$, i.e. $k=1,2,4$, the number $2p^k$ can be written in the form   $2p^k=q^{2^l}+1, l>0$. Moreover, in the case $k=4$, we can state the result analog to Theorem~\ref{tm2.2}:


\begin{theorem}\label{tm2.2_s4}
Let $t>0$ and $\{1,b,c\}$ be $D(-1)$-triple in the ring $\Z[\sqrt{-t}]$.
If $b=2 p^4$, where $p$ is an odd prime, then $c\in \Z$.
\end{theorem}

The proof of Theorem~\ref{tm2.2_s4} follows the same steps as the proof of Theorem~\ref{tm2.2}(ii), (iii), so we will omit it.

\begin{remark}\label{rem_21}
The statement of Remark 2 is valid in the case of Theorem \ref{tm2.2_s4}, too.
\end{remark}

In proving results of this section we will use the following result of Filipin, Fujita
and Mignotte from \cite{FiFuMi} on $D(-1)$-quadruples in integers.

\begin{lemma} [\text{\cite[Corollary 1.3]{FiFuMi}}]\label{lema_FiFuMi}
Let $r$ be a positive integer and let $b=r^2+1$. Assume
that one of the following holds for any odd prime $p$
and a positive integer $k$:
\be
b=p, \quad  b=2p^k, \quad  r=p^k, \quad  r= 2p^k.
\ee
Then the system of Diophantine equations
\be
y^2-bx^2&=&r^2,\\
z^2-cx^2&=&s^2
\ee
has only the trivial solutions $(x,y,z)=(0,\pm r,\pm s)$,
where $s$ is such that $(t,s)$ is  a positive solution of $t^2-b s^2=r^2$
and $c=s^2+1$. Furthermore, the $D(-1)$-pair $\{1,b\}$ cannot
be extended to a $D(-1)$-quadruple.
\end{lemma}

First we prove the following result.

\begin{theorem}\label{tm(ii)_1}
If $p$ is an odd prime and $t \equiv 0 \pmod 2$, then there does not exists a $D(-1)$-quadruple
of the form $\{1,2p^k,c,d\}$ in $\Z[\sqrt{-t}]$.
\end{theorem}

\noindent
\pf Let $t \equiv 0 \pmod 2$. We have that $t\nmid 2p^k-1$. Therefore, if we suppose
that $\{1,2p^k,c,d\}$ is a $D(-1)$-quadruple in $\Z[\sqrt{-t}]$, then according to Remark~\ref{rem_2} and \ref{rem_21} we obtain $c,d \in \N$.
This means that there exist integers $x_1,y_1,u_1,v_1,w_1$, such that
\be
c-1=x_1^2, d-1=y_1^2, 2p^kc-1=u_1^2, 2p^kd-1=v_1^2, cd-1=w_1^2,
\ee
or at least one of  $c-1, d-1, 2p^kc-1, 2p^kd-1, cd-1$ is equal to $-tw_2^2$, for an integer $w_2$.

The first possibility leads to contradiction with Lemma~\ref{lema_FiFuMi}, i.e.,
a $D(-1)$-pair $\{1,2p^k\}$, cannot be extended to a $D(-1)$-quadruple in integers,
while the second one contradicts to $c,d \in \N$.

\hfill $\Box$

In what follows, our main goal is to obtain some results for odd $t$'s. Thus, let us consider the case of $t \equiv 1 \pmod 2$. We have the following result:

\begin{theorem}\label{tm(ii)_2_novi}
Let $2p^k=q^{2^l}+1, l>0$, where $p$ and $q$ are odd primes.
\begin{itemize}
\item[$($i$)$] If $t\in\{1,q^2, \dots, q^{2^l-2}, q^{2^l}\}$, then there exist infinitely many $D(-1)$-quad\-ruples of the form
$\{1,2p^k,-c,d\}$, $c,d>0$ in $\Z[\sqrt{-t}]$.

\item[$($ii$)$] If $t\in\{q,q^3, \dots,q^{2^l-3}, q^{2^l-1}\}$, then there does not exists a $D(-1)$-quadruple
of the form $\{1,2p^k,c,d\}$ in $\Z[\sqrt{-t}]$.
\end{itemize}
\end{theorem}

Before we start with proving Theorem~\ref{tm(ii)_2_novi}, we recall the following result.

\begin{lemma}[\textnormal{\cite[Lemma 3]{DujeCambridge}}] \label{lemma2_cet}
If $\{a,b,c\}$ is a Diophantine triple with the property $D(l)$ and $ab+l=r^2$, $ac+l=s^2$, $bc+l=t^2$,
then there exist integers $e,x,y,z$ such that
\[
ae+l^2=x^2, be+l^2=y^2, ce+l^2=z^2
\]
and
\[
c=a+b+\frac{e}{l}+\frac{2}{l^2}(abe+rxy).
\]
\end{lemma}
Moreover, $e=l(a+b+c)+2abc-2rst$, $x=at-rs$, $y=bs-rt$, $z=cr-st$.

To prove the next proposition, which will be used in proving Theorem~\ref{tm(ii)_2_novi}, we will use Lemma~\ref{lemma2_cet} for $l=-1$.

\begin{proposition}\label{prop2_6} Let $m,n>0$ and $b=n^{2}+1$.
If $m|n$ and $t\in\{1,m^2, n^{2}\}$, then there exist infinitely many $D(-1)$-quadruples of the form
$\{1,b,-c,d\}$, $c,d>0$ in $\Z[\sqrt{-t}]$.
\end{proposition}

\noindent
{\bf Proof}:
Since $\Z[n i]$ is a subring of $\Z[m i]$,
it suffices to prove the statement for $t=n^{2}$.
Thus, suppose that there exist $x,y\in \Z$ such that
\be
-c-1&=&-n^{2} x^2=(n xi)^2,\\
-bc-1&=&-n^{2} y^2=(n yi)^2.
\ee
Eliminating $c$, we obtain Pellian equation
\bee
y^2-(n^{2}+1)x^2=-1. \label{eq_Z[2i]}
\eee
All positive solutions of the equation \eqref{eq_Z[2i]} are given by
\be
x=x_j &=& \frac{\sqrt{n^{2}+1}}{2 (n^{2}+1)}\left((n +\sqrt{n^{2}+1})^{2j-1} - (n -\sqrt{n^{2}+1})^{2j-1}\right),\\
y=y_j &=& \frac{1}{2}\left((n +\sqrt{n^{2}+1})^{2j-1} + (n -\sqrt{n^{2}+1})^{2j-1}\right), \quad j \in \N.
\ee
Therefore, for any $j \in \N$ and
$
c=c_j=n^{2} x_j^2-1,
$
the set  $\{1,b,-c\}$ is a $D(-1)$-triple in $\Z[\sqrt{-t}]$. If we apply Lemma~\ref{lemma2_cet} on
that triple, we obtain positive integers
\be
d_{+,-}&=&\pm 2 n^{3} x_j y_j +(2n^{2}+1)c + n^{2} + 2,
\ee
such that
\be
d_{+,-}-1&=&\left(n^{2} x_j \pm n y_j\right)^2,\\
b d_{+,-}-1 &=&\left(n(n^{2}+1) x_j \pm n^{2} y_j\right)^2,\\
-c d_{+,-}-1&=&\left(n c i \pm n^{2} x_j y_j i\right)^2.
\ee
Thus the sets $\{1,b,-c,d_+\}, \{1,b,-c,d_-\}$ are $D(-1)$-quadruples in $\Z[\sqrt{-t}]$.

\hfill $\Box$

Now, we are able to prove Theorem~\ref{tm(ii)_2_novi}.

\medskip
\noindent
{\bf Proof of Theorem~\ref{tm(ii)_2_novi}:}

\noindent
Let $l\geq 0$.

i) Suppose that $t \in \{1, q^2, \dots, q^{2^l-2}, q^{2^l}\}$.
By Proposition~\ref{prop2_6} there exists infinitely many $D(-1)$-quadruples of the form
$\{1,2p^k,-c,d\}$, $c,d>0$ in $\Z[\sqrt{-t}]$.

ii) Let us assume that $t\in\{q,q^3, \dots,q^{2^l-3}, q^{2^l-1}\}$.
In this case, the equation \eqref{eq7} is equivalent to
\bee\label{eq8}
x^2-(q^{2^l}+1) y^2=-q^s,
\eee
where $s$ is an odd integer and $0<s\leq 2^l-1$. Theorem~\ref{tm1} implies that the equation \eqref{eq8} has no integer solutions. Therefore, if  $\{1,2p^k,c,d\}$ is $D(-1)$-quadruple in $\Z[\sqrt{-t}]$,
then $c,d>0$. By the same argumentation as in Theorem~\ref{tm(ii)_1} we conclude that such quadruple does not exist.

\hfill $\Box$

\bigskip
\noindent
Andrej Dujella, Department of Mathematics, University of Zagreb, Bijeni\v{c}ka cesta 30, 10000 Zagreb, Croatia\\
e-mail: duje@math.hr

\vspace*{5mm}
\noindent
Mirela Juki\'c Bokun, Department of Mathematics, University of Osijek, Trg Ljudevita Gaja 6, 31000 Osijek, Croatia\\
e-mail: mirela@mathos.hr

\vspace*{5mm}
\noindent
Ivan Soldo, Department of Mathematics, University of Osijek, Trg Ljudevita Gaja 6, 31000 Osijek, Croatia\\
e-mail: isoldo@mathos.hr

\end{document}